\documentclass[12pt]{article}
\usepackage{hyperref}
\usepackage{dsfont}
\usepackage{amssymb}
\newcommand{\NN}{\mathds{N}}
\newcommand{\ZZ}{\mathds{Z}}
\newcommand{\QQ}{\mathds{Q}}
\newcommand{\RR}{\mathds{R}}
\begin{document}
\begin{center}
{\bf One-Dimensional Packing: \\
Maximality Implies Rationality} \\
\ \\
James Propp \\
October 30, 2017
\end{center}

\begin{abstract}
\noindent
{\sc Abstract:} Every set of natural numbers
determines a generating function convergent for $q \in (-1,1)$
whose behavior as $q \rightarrow 1^-$ determines a germ.
These germs admit a natural partial ordering
that can be used to compare sizes of sets of natural numbers
in a manner that generalizes both cardinality of finite sets
and density of infinite sets.
For any finite set $D$ of positive integers, call a set $S$ ``$D$-avoiding'' 
if no two elements of $S$ differ by an element of $D$.
It is shown that any $D$-avoiding set 
that is maximal in the class of $D$-avoiding sets (with respect to germ-ordering)
is eventually periodic.
This implies an analogous result for packings in $\NN$.
It is conjectured that for all finite $D$ 
there is a unique maximal $D$-avoiding set.
\end{abstract}

\bigskip

\noindent
{\bf 1. Introduction}

\bigskip

This article is concerned with two related kinds 
of optimization problems in $\NN$:
packing problems and distance-avoidance problems.
In the former, we are given a nonempty set $B \subseteq \NN = \{0,1,2,\dots\}$
and we wish to find a collection of disjoint translates of $B$
whose union is as big a subset of $\NN$ as possible.
In the latter, we are given a finite set $D$ of positive integers
and we wish to find as big a set $S \subseteq \NN$ as possible
such that no two elements of $S$ differ by an element of $D$.
In both cases, the crucial issue is defining what ``as big as possible'' should mean.

For instance, consider the distance-avoidance problem with $D = \{3,5\}$.
Three $D$-avoiding sets are $S_0 = \{0,2,4,6,\dots\}$,
$S_1 = \{1,3,5,7,\dots\}$, and $S_2 = \{0,1,2,8,9,10,\dots\}$
(note that the third set is obtained via an obvious general algorithm 
for greedily constructing $D$-avoiding sets for arbitrary $D$).
In terms of subset-inclusion, all three sets are maximal:
none of them can be augmented
without violating the $D$-avoidance property.
We will say $S_0$ is ``bigger'' than $S_1$, which is in turn ``bigger'' than $S_2$,
in the sense that
$\sum_{n \in S_0} q^n > \sum_{n \in S_1} q^n > \sum_{n \in S_2} q^n$
for all $q<1$ sufficiently close to 1.
That is, we propose to measure of the size of a set $S \subseteq \NN$
by forming the generating function $S_q := \sum_{n \in S} q^n$
and examining its germ ``at $1^-$''.

For example:
\begin{enumerate}
\item
If $S$ is finite, $S_q = |S| + o(1)$, or equivalently,
$S_q \rightarrow |S|$ as $q \rightarrow 1^-$;
if $S$ is infinite, $S_q$ diverges as $q \rightarrow 1^-$.
\item
If $S$ is infinite with density $\alpha$,
$S_q = \alpha \frac{1}{1-q} + o(\frac{1}{1-q})$.
\item
In particular, if $S = \{a,a+d,a+2d,\dots\}$
with $a \geq 0$ and $d > 0$,
$S_q = (\frac{1}{d}) \frac{1}{1-q} + (\frac{d-1-2a}{2d}) + O(1-q)$.
\end{enumerate}

This approach is related to Abel's method of evaluating divergent series;
its application to measuring sets of natural numbers is (apparently) new,
but it is likely to hold little novelty for analytic number theorists,
who have long used the philosophically similar
but technically more recondite notion of Dirichlet density
to measure sets of primes.
Our definition also has thematic links to work from 
the earliest days in the study of infinite series.
For instance, Grandi's formula $1-1+1-1+1-1+\dots=\frac12$ corresponds to the fact 
that the germ of $(2\NN)_q$ exceeds the germ of $(2\NN+1)_q$ by $\frac12 + O(1-q)$.
while Callet's formula $1+0-1+1+0-1+\dots=\frac23$ corresponds to the fact
that the germ of $(3\NN)_q$ exceeds the germ of $(3\NN+2)_q$ by $\frac23 + O(1-q)$.

Our approach resembles the sort of ``tame nonstandard analysis''
in which $\RR$ is replaced by the ordered ring $\RR(x)$ 
where $1/x$ is a formal infinitesimal
(also known as ``the ring of rational functions ordered at infinity'');
our ordering of rational functions corresponds to that of $\RR(x)$
if one identifies $1/x$ with $1-q$.

Theorems 2 and 4 show that 
for both packing problems and distance-avoidance problems in $\NN$,
every optimal (that is, germ-maximal) solution is eventually periodic.
The proof we give may seem surprisingly complicated, given that
the corresponding periodicity property for maximum-density packings
and maximum-density distance-avoiding sets is fairly easy.
This discrepancy is explained by the fact that
the germ-topology does not admit compactness arguments.

We conjecture that for both the packing and distance-avoidance problems,
there is a unique optimum subset of $\NN$ (guaranteed to be eventually periodic).

The motivation for this work was the study of disk packings.
It is our hope that the approach taken here will ultimately lead to results
establishing a strong kind of uniqueness
for optimal sphere-packings in dimensions 2, 8, and 24.
(See [Co] for a survey of the recent breakthroughs
in the study of 8- and 24-dimensional sphere-packing.)
We also hope that the germ approach will have relevance
to the study of densest packings in other dimensions.

For other approaches to measuring efficiency of packings, see [Ku].
The most sophisticated of these approaches is that of Bowen and Radin [Bo];
their ergodic theory approach has attractive features
(for instance, it works in spaces with nonamenable symmetry groups),
but it does not seem to work so well when the region being packed
is not the entire space.  Packings in $\NN$
could be viewed as special packings of $\RR^{\geq 0}$;
the lack of symmetry makes it hard to apply the constructions of Bowen and Radin.

See also [Be], [Bl], [Ch], and [Ka] for work on measuring sizes of sets
bearing some philosophical similar to ours.

\bigskip

\noindent
{\bf 2. Statement of main theorem}

\bigskip

Recall that a subset $S$ of $\NN = \{0,1,2,\dots\}$
is eventually periodic iff there exist $N \in \NN$ and $d \geq 1$
such that for all $n \geq N$, $n \in S$ iff $n+d \in S$.
It is easy to show that $S$ is eventually periodic
if and only if its {\bf generating function} $S_q := \sum_{n \in S} q^n$
is a rational function of $q$.  We call such sets $S$ {\bf rational}.
(Note that this usage coincides with the notion of rationality
for subsets of a monoid in automata theory, specialized to the monoid $\NN$.)
If $S$ is a finite set, then $S$ is rational and $S_q$ is a polynomial.
If $S$ is rational and infinite, then $S_q$ has a simple pole at 1,
and letting $t = 1-q$ 
we can expand $S_q$ as a Laurent series $\sum_{n \geq -1} a_n t^n$
where $a_{-1}$ is the density of $S$.
This series converges for all $q$ in $(-1,1)$, 
though we will only care about $q$ in $(0,1)$.

Given two sets of natural numbers $S$ and $S'$ (not necessarily rational),
write $S \preceq S'$ iff there exists $\epsilon > 0$ such that
$S_q \leq S'_q$ for all $q$ in the interval $(1-\epsilon,1)$;
we say that $S'$ {\bf dominates} $S$ in the germ-ordering.
The partial ordering $\preceq$ (which we call the {\bf germ-ordering at $1^-$})
is a total ordering on the rational subsets of $\NN$
that refines the preorder given by comparing density.
Also, if two sets have finite symmetric difference they are $\preceq$-comparable.
(Both of these assertions are consequences of the fact that
the sign of a polynomial can oscillate only finitely many times.)
In the case where $S$ and $S'$ are finite,
the germ-ordering refines ordering by cardinality;
when the finite sets $S$ and $S'$ have the same cardinality $n$,
the germ-ordering refines lexicographic ordering of subsets of $\NN$ of size $n$.
(When $S,S'$ are eventually periodic infinite sets of the same density $c$,
there is also a combinatorial criterion for deciding which of $S,S'$ is larger,
though it is more complicated.)

The germ-ordering has the ``outpacing property'' [Ka]:
if for all sufficiently large $k$
the $k$th element of $S$ is less than or equal to the $k$th element of $S'$,
then $S \succeq S'$.

We mention that, although $\preceq$ is a total ordering for rational subsets of $\NN$, 
the same is not true for unrestricted subsets of $\NN$;
for instance, if $S$ is the set of natural numbers 
whose base ten expansion has an even number of digits
and $S'$ is its complement,
then it can be shown that $S$ and $S'$ are $\preceq$-incomparable.

Given a finite nonempty subset $B$ of $\NN$ (a {\bf packing body}),
say that a set $T \subset \NN$ is a {\bf translation set for $B$}
iff the translates $B+n$ ($n \in T$) are disjoint. 
If $T$ is a translation set,
the generating function of $\cup_{n \in T} (B+n)$
is just the product of the generating function of $T$
and the generating function of $B$;
so if $T$ and $T'$ are translation sets,
$T \preceq T'$ iff $\cup_{n \in T} (B+n) \preceq \cup_{n \in T'} (B+n)$.

\medskip

\noindent
{\bf Conjecture 1:} For every packing body $B$, 
there is a unique germ-maximal translation set for $B$, and it is rational.
That is, there is a translation set $T^*$ such that $T^*$ is rational
and such that $T \preceq T^*$ for every translation set $T$.

\medskip

\noindent
This Conjecture is easy to prove for many specific packing bodies,
such as $\{0,1,\dots\,k-1\}$ for arbitrary $k$
(see Section 4), but we do not have a general proof.
Theorem 2 is the best result we currently have
that applies to all packing bodies $B$.

\medskip

\noindent
{\bf Theorem 2:} For every packing body $B$,
every germ-maximal translation set is rational.
That is, if $T^*$ is a translation set with the property that
there exists no translation set $T \succ T^*$, then $T^*$ is rational.

\medskip

\noindent
We hope to (but cannot yet) prove
that the collection of translation sets for $B$ contains a maximal element;
it is a priori conceivable that there exist translation sets
$T_1 \prec T_2 \prec T_3 \prec \dots$
but no translation set that dominates them all.
Thus Theorem 2 does not immediately imply Conjecture 1.

Our proof of Theorem 2 goes by way of a shift of context
from the packing problem to the forbidden distance problem
(which in $\NN$ might with equal aptness be called the forbidden difference problem).
The condition that $T$ is a translation set for $B$
is equivalent to the condition that
the difference set $T-T = \{x-y: x,y \in T\}$ has no element in common with 
the difference set $B-B = \{x-y: x,y \in B\}$ other than 0.
Thus the problem of finding the germ-maximal translation set for the packing body $B$ 
is a special case of the problem of finding the germ-maximal set $T \subseteq \NN$
that has no differences in the finite set $D_B$
where $D_B$ is the set of positive elements of $B-B$.
More generally, for any finite set $D$ of positive integers,
say that $S \subseteq \NN$ is {\bf $D$-avoiding}
if there exist no two elements in $S$ that differ by an element of $D$.
In this setting we can broaden Conjecture 1 and Theorem 2.

\medskip

\noindent
{\bf Conjecture 3:} For every finite set $D$ of positive integers,
there is a unique germ-maximal $D$-avoiding set $S^*$ and it is rational.

\medskip

\noindent
{\bf Theorem 4:} For every finite set $D$ of positive integers,
every germ-maximal $D$-avoiding set is rational.

\medskip

\noindent
Of course Conjecture 3 implies Conjecture 1 and Theorem 4 implies Theorem 2.

The conclusions of Theorem 2 and Theorem 4 cannot be strengthened to assert that
the set must be periodic, as is demonstrated by the following example
(jointly found with Aaron Abrams, Henry Landau, Zeph Landau,
Jamie Pommersheim, and Alexander Russell):
Let $B = \{0,4,11\}$ and $D = \{4,7,11\}$ (the set of positive elements of $B-B$).
The germ-maximal periodic $D$-free subset of $\NN$
is the period-3 set $\{0, 3, 6, 9, 12, 15, 18, \dots\}$
but the eventually periodic set $\{0, 1, 3, 6, 9, 15, 18, \dots\}$
(in which 12 is replaced by 1) is infinitesimally larger,
and is indeed the germ-maximal $D$-free subset of $\NN$.
It follows that the germ-densest packing of $B$
is eventually periodic but not periodic.
(Details will appear elsewhere.)
Note that this example undermines Conjectures 9 and 10 from
the earlier posted version of this paper.
It is possible that every one-dimensional packing problem
has a periodic solution that is optimal modulo infinitesimals
(that is, up to germs that are $o(1)$ as $q \rightarrow 1$).
Abrams et al.\ also showed that the set $2\NN$
is the germ-maximal $D$-avoiding set for $D=\{3,5\}$.

\bigskip

\noindent
{\bf 3. Proof of main theorem}

\bigskip

Our approach to proving Theorem 4 uses a block coding 
of the kind often employed in dynamical systems theory.
Let $m = \max(D) + 1$ and replace the indicator sequence of $S$
(an element of $\{0,1\}^\NN$) by a symbolic sequence
using a block code of block length $m$,
with an alphabet containing (at most) $2^m$ symbols, which we will call {\bf letters}.
More concretely, if the indicator sequence of $S$ is written as $(b_0,b_1,b_2,\dots)$
(where $b_n$ is 1 or 0 according to whether $n \in S$ or $n \not\in S$),
then we define the {\bf $m$-block encoding} of $(b_0,b_1,b_2,\dots)$
to be $(w_0,w_1,w_2,\dots)$
where the letter $w_n$ is the $m$-tuple $(b_n,b_{n+1},\dots,b_{n+m-1})$;
we call $w_n$ a {\bf consonant} or a {\bf vowel}
according to whether $b_n = 1$ or $b_n = 0$
(conditions that align with the respective cases $n \in S$ and $n \not \in S$).
Say that a letter $\alpha=(b_1,\dots,b_m)$ in $\{0,1\}^m$ is {\bf legal}
if the set $\{i: b_i=1\}$ is $D$-avoiding;
we let $\cal{A}$ be the set of legal letters.
Given two letters $\alpha$ and $\alpha'$ in $\cal{A}$,
say that $\alpha'=(b'_1,\dots,b'_m)$ is a {\bf successor} of $\alpha=(b_1,\dots,b_m)$
iff $b'_i = b_{i+1}$ for $1 \leq i \leq m-1$.
For every set $S \subseteq \NN$,
the associated block-encoding $w = (w_0,w_1,w_2,\dots)$
has the property that for all $n \geq 1$,
$w_n$ is a successor of $w_{n-1}$;
$S$ is $D$-avoiding if and only if $w$
has the additional property that every letter $w_n$ is legal.
Call such an infinite word $(w_0,w_1,w_2,\dots)$ {\bf $D$-legal}.
Finding a germ-maximal $D$-avoiding set is equivalent to finding 
a $D$-legal infinite word 
for which the set of locations of consonants is germ-maximal.
We write $w \preceq w'$ iff the associated sets $S,S'$ satisfy $S \preceq S'$.

Suppose $S$ is some $D$-avoiding subset of $\NN$
that is germ-maximal in the collection of $D$-avoiding subsets of $\NN$.
Let $w=(w_0,w_1,\dots)$ be the associated infinite word in $\cal{A}^\NN$.
Assume for simplicity that the letter $w_0 = \alpha$ occurs infinitely often in $w$.
(The last paragraph of the proof addresses what happens
if this assumption fails.)

Let $K = \{k \in \NN: w_k = \alpha \} = \{k_0, k_1, k_2, \dots\}$,
where $k_0 = 0$ and $k_0 < k_1 < k_2 < \dots$.
This divides up the infinite word $w$ into infinitely many subwords
$(w_{k_0},w_{k_0+1},\dots,w_{k_1-1})$,
$(w_{k_1},w_{k_1+1},\dots,w_{k_2-1})$,
$(w_{k_2},w_{k_2+1},\dots,w_{k_3-1})$, \dots.
Each of these finite words $(w_{k_i},w_{k_i+1},\dots,w_{k_{i+1}-1})$
is associated with the word
$c_k := (w_{k_{i-1}},w_{k_{i-1}+1},\dots,w_{k_{i}-1},w_{k_{i}})$
that both begins and ends with the letter $\alpha$;
define a {\bf circular word} as a word whose first and last letters are the same.
(Note that we are not modding out by cyclic shift of such words.)
Let $\cal{C}$ be the set of all circular words beginning and ending with $\alpha$.
We define the {\bf length} of a circular word to be the number of letters it contains,
counting its first and last letter as a single letter.
(Thus, if $\alpha$, $\beta$, and $\gamma$ are letters, 
the circular word $\alpha \beta \gamma \alpha$ has length 3.)
If $c \in \cal{C}$ has length $a$ and $c' \in \cal{C}$ has length $a'$,
let $c\!:\!c'$ denote the circular word of length $a+a'$ in $\cal{C}$
obtained by concatenating $c$ and $c'$
(where the final $\alpha$ in $c$ gets identified with the initial $\alpha$ in $c'$).
The operation $:$ is associative, and indeed,
the word $w$ itself can be written as $c_1\!:\!c_2\!:\!c_3\!:\!\dots$,
where the circular words $c_i$ are {\bf primitive}
(i.e., each $c_i$ contains $\alpha$ only at the beginning and at the end).
We also use ``:'' to denote concatenation of noncircular words. 

Every circular word $c \in \cal{C}$ is associated 
with a polynomial $P_c = P_c(q)$
(sometimes we will omit the subscript or will write $P_i$ to mean $P_{c_i}$)
whose degree is at most the length $a$ of the circular word $c$
and whose coefficients are 0's and 1's according to whether 
the respective letters in the circular word are vowels or consonants;
we call $P_c$ the {\bf generating function} of $c$.
So if $w = c_1 \!:\! c_2 \!:\! c_3 \!:\! \dots$ is the $D$-legal infinite word
representing the $D$-avoiding set $S$, $S_q$ can be written as 
$P_1 + q^{a_1} P_2 + q^{a_1+a_2} P_3 + \dots = P_1 + A_1 P_2 + A_1 A_2 P_3 + \dots$
where $a_i$ is the length of $c_i$ and $A_i$ is $q^{a_i}$.

For any circular word $c$ with length $a$, we define $|c| := P_c(q) / (1 - q^a)$;
it is equal to the generating function 
of the infinite periodic word $c\!:\!c\!:\!c\!:\!c\!:\!\dots$.
Given two periodic words $c,c'$ in $\cal{C}$ (possibly of different lengths),
write $c \preceq c'$ iff $|c| \preceq |c'|$; 
call this the germ-ordering on circular words.
We have $|c| = |c'|$ iff $c\!:\!c\!:\!c\!:\!\dots = c'\!:\!c'\!:\!c'\!:\!\dots$.

The following two lemmas are the linchpins of the proof of Theorem 4.

\medskip

\noindent
{\bf Lemma 5:} If $c \preceq c'$, then
$c \preceq c\!:\!c' \preceq c'\!:\!c \preceq c'$.

\medskip
\noindent
{\bf Proof:}
Write $|c| = P/(1-A)$ and $|c'| = P'/(1-A')$;
we also have $|c\!:\!c'| = (P + AP') / (1-AA')$
and $|c'\!:\!c| = (P' + A'P) / (1-AA')$.
The stipulated relation $c \preceq c'$ is equivalent to
$P/(1-A) \preceq P'/(1-A')$, or
\begin{equation}
P(1-A') \preceq P'(1-A);
\end{equation}
the desired relations $c \preceq c\!:\!c'$, 
$c\!:\!c' \preceq c'\!:\!c$, and $c'\!:\!c \preceq c'$
are respectively equivalent to
\begin{equation}
P/(1-A) \preceq (P + AP') / (1-AA'),
\end{equation}
\begin{equation}
(P + AP') / (1-AA') \preceq (P' + A'P) / (1-AA'), \ \mbox{and}
\end{equation}
\begin{equation}
(P' + A'P) / (1-AA') \preceq P'/(1-A').
\end{equation}
To prove (2), note that (by cross-multiplying, expanding, and cancelling terms)
we can write it equivalently as $- AA'P \preceq AP' - AP - AAP'$,
which is just (1) multiplied by $A$.
The two denominators in (3) are identical,
so (3) is equivalent to $P + AP' \preceq P' + A'P$,
which in turn is equivalent to (1).
The proof of (4) is similar to the proof of (2).
$\blacksquare$

\medskip

Note that the proof also tells us that if $c \prec c'$, 
then $c \prec c\!:\!c' \prec c'\!:\!c \prec c'$.

\medskip

\noindent
{\bf Lemma 6:} If the concatenation $w = c_1\!:\!c_2\!:\!c_3\!:\! \dots$ 
is germ-maximal in the set of $D$-legal words,
then we must have $c_1 \succeq c_2 \succeq c_3 \succeq \dots$ in the germ-ordering.

\medskip

\noindent
{\bf Proof:}
We will show that $c_1 \succeq c_2$ 
since that contains the idea of the general argument.
If $c_1 = c_2$ there is nothing to prove, so assume $c_1 \neq c_2$,
and let $w' = c_2\!:\!c_1\!:\!c_3\!:\! \dots$, which must be $D$-legal if $w$ is
(indeed, the whole reason for the block coding was to make this claim true).
The sets $S$ and $S'$ respectively associated with $w$ and $w'$
have finite symmetric difference, so $w$ and $w'$ must be comparable.
Since we are assuming $w$ is germ-maximal,
we must have $w \succeq w'$ in the germ ordering.
That is, we must have $$P_1 + A_1 P_2 \succeq P_2 + A_2 P_1$$
(all the later terms match up and cancel).
But this is equivalent to $P_1 / (1 - A_1) \succeq P_2 / (1 - A_2)$, 
so $c_1 \succeq c_2$ as claimed. $\blacksquare$

\medskip

\noindent
{\bf Proof of Theorem 4:}
By an easy pigeonhole argument, for all $N$ there must exist $i,j \geq N$ with $i < j$
such that the sum of the lengths of the words $c_i$, $c_{i+1}, \dots, c_j$
is a multiple of the length of $c_1$, say $r$ times the length of $c_1$.
Let $w'$ be the word obtained from $w$
by replacing the $j-i+1$ letters $c_i$, $c_{i+1}, \dots, c_j$
by $r$ occurrences of the letter $c_1$.
Let $S$ and $S'$ be the sets associated with $w$ and $w'$, respectively.
Lemma 6 tells us that $c_1 \succeq c_i \succeq c_{i+1} \succeq \dots \succeq c_j$,
so repeated application of Lemma 5 gives
$|c_1\!:\!c_1\!:\!\dots\!:\!c_1| \succeq |c_{i}\!:\!c_{i+1}\!:\!\dots\!:\!c_{j}|$.
If strict inequality holds, then $w' \succ w$, contradicting maximality of $w$.
(Here we use the fact that the difference $S'_q - S_q$
can be expressed as $1-q^n$ times
$|c_1\!:\!c_1\!:\!\dots\!:\!c_1| - |c_{i}\!:\!c_{i+1}\!:\!\dots\!:\!c_{j}|$,
where $n$ is the common value of $ra_1$ and $a_i+a_{i+1}+\dots+a_j$.)
So we must have
$|c_1\!:\!c_1\!:\!\dots\!:\!c_1| = |c_{i}\!:\!c_{i+1}\!:\!\dots\!:\!c_{j}|$,
implying that $c_i,c_{i+1},\dots,c_j$ are all the circular word $c_1$.
Since the circular words $c_i$ are in germ-decreasing order,
this means that $c_1,c_2,\dots,c_N$ are all equal.
Since this is true for all $N$, we must have $w = c_1\!:\!c_1\!:\!\dots c_1$;
that is, $w$ is periodic.

The above argument was predicated on the assumption that $\alpha$
occurs infinitely often.
If this assumption fails, then a version of the argument still goes through,
but it is slightly more complicated;
one finds the smallest $i$ for which
the letter $w_i$ occurs infinitely often in $w$
(guaranteed to exist),
and then one applies the preceding argument
to the letters $w_i,w_{i+1},w_{i+2},\dots$,
ignoring the letters $w_0,\dots,w_{i-1}$.
Instead of concluding that $w$ is periodic,
we obtain the weaker conclusion that $w$ is eventually periodic.
$\blacksquare$

\bigskip
\noindent
{\bf 4. Existence and uniqueness in special cases}

\bigskip

In the case where $D = \{1,2,\dots,k-1\}$ for some $k \geq 1$, 
it is easy to give a direct proof of existence and uniqueness
of a maximal $D$-avoiding set, namely $S^* = \{0,k,2k,\dots\}$.
$S^*$ dominates every $D$-avoiding set $S$ in the sense that, 
writing $S = \{s_1,s_2,s_3,\dots\}$ with $s_1 < s_2 < s_3 < \dots$,
we have $s_1 \geq 0$, $s_2 \geq k$, $s_3 \geq 2k$, etc.

The result for $D = \{1,2,\dots,k-1\}$ is also implied by a more general result:

\medskip

\noindent
{\bf Theorem 7:} Suppose that for every letter $\alpha \in \cal{A}$
there exists a circular word $c_\alpha^*$ whose first letter is $\alpha$, such that
$c_\alpha^* \succeq c$ for every circular word $c$ whose first letter is $\alpha$.
Then every $D$-legal word $w$ is dominated by
an eventually periodic word whose repetend is one of the circular words $c_\alpha$
and whose ``pre-repetend'' does not contain any repeated letters.

\medskip
\noindent
{\bf Proof:}
Let $\alpha$ be some letter that occurs in $w$ infinitely often,
and write $w$ as $b \!:\! c_1 \!:\! c_2 \!:\! c_3 \!:\! \dots$,
where $c_1,c_2,c_3,\dots$ all start with $\alpha$.
Let $n$ be the length of $c_\alpha^*$.
Consider the lengths of the truncated words
$b \!:\! c_1 \!:\! c_2 \!:\! c_3 \!:\! \dots \!:\! c_m$ 
(for $m \geq 1$) mod $n$;
some residue class must be represented infinitely often,
so we can find $i_1 < i_2 < \dots$ such that for all $k$,
the length of $c_{i_k} \!:\! c_{i_k+1} \!:\! \dots \!:\! c_{i_{k+1}-1}$
is a multiple of $n$.
Then we can replace each such stretch of $w$
by a succession of occurrences of $c_\alpha^*$ of the exact same total length, satisfying 
$c_\alpha^* \!:\! c_\alpha^* \!:\! \dots \!:\! c_\alpha^* \succeq
c_{i_k} \!:\! c_{i_k+1} \!:\! \dots \!:\! c_{i_{k+1}-1}$;
this results in an eventually periodic word $w^*$
whose repetend is $c_\alpha^*$, satisfying $w^* \succeq w$.

Now we must show that $w^*$ is in turn dominated
by a word whose pre-repetend does not contain any repeated letters.
Suppose the pre-repetend contains two occurrences of the letter $\alpha'$.
Write $w^*$ as $d \!:\! e \!:\! f \!:\! c_\alpha^* \!:\! c_\alpha^* \!:\! \dots$,
where the finite words $e$ and $f$ both start with $\alpha'$. 
Let $A$, $B$, $C$, and $D$ be $q$ to the power of
the lengths of $c_\alpha^*$, $d$, $e$, and $f$, respectively,
and let $P$, $Q$, $R$, and $S$ be the polynomial generating functions
of $c_\alpha^*$, $d$, $e$, and $f$, respectively.
Since rational functions in $q$ form a totally ordered set under germ-ordering,
we must have $R \preceq (1-C)(S+DP/(1-A))$ or $R \succeq (1-C)(S+DP/(1-A))$ or both.
In the former case, we have $Q + BS + BDP/(1-A) \succeq Q + BR + BCS + BCDP/(1-A)$,
so that $d \!:\! f \!:\! c_\alpha^* \!:\! c_\alpha^* \!:\! \dots$ dominates $w^*$
(that is, we make the word $w^*$ bigger by removing the subword $e$);
in the latter case, we have
$Q + BR/(1-C) \succeq Q + BR + BCS + BCDP/(1-A)$,
so that $d \!:\! e \!:\! e \!:\! e \!:\! \dots$ dominates $w^*$
(that is, we make $w^*$ bigger by putting in infinitely many $e$'s).
Either way, we get an infinite word with strictly shorter pre-repetend,
so if we iterate the procedure as needed,
we must eventually arrive at an infinite word
whose pre-repetend contains no repeated letters. 
$\blacksquare$

\bigskip


We mentioned in the introduction that germs do not come with a nice topology. 
As an illustration of this (related to the famous Ross-Littlewood Paradox),
consider the sequence of sets $S_n = \{n,n+1,\dots,10n\}$;
we have $S_1 \prec S_2 \prec S_3 \prec \dots$,
but it is unclear what the limit of the $S_n$'s should be.
Surely it is not the pointwise limit of the sets, since that is the null set!
One way to understand what is going on here
is to note that, even though for each $n$ there exists $\epsilon_n > 0$
such that $(S_n)_q < (S_{n+1})_q$ for all $q$ in $(1-\epsilon_n,1)$,
we have $\inf \epsilon_n = 0$,
so that the intersection of the intervals $(1-\epsilon_n,1)$ is empty.

This sort of situation comes into play when one tries 
to prove Conjecture 3 by showing that
$c \succeq c_1,c_2,c_3,\dots$ implies
$c\!:\!c\!:\!c\!:\!\dots \succeq c_1\!:\!c_2\!:\!c_3:\dots$.
If we take $\epsilon_n$ satisfying $|c| \geq |c_n|$ for all $q$ in $(1-\epsilon_n,1)$,
and the infimum of the $\epsilon_n$ not known to be positive,
then the obvious approach to proving the implication fails.

\bigskip

\noindent
{\bf 5. Truncating the germs}

\bigskip

In our approach, a rational set $S \subseteq \NN$ is replaced by the power series
$\sum_{n \in S} q^n$, which is rewritten as the Laurent series $\sum_{n \geq -1} a_n(1-q)^n$,
and the coefficients $a_{-1},a_0,a_1,a_2,\dots$ are used to put
a total ordering on the rational sets.
The coefficients $a_n$ carry finer and finer information as $n$ increases,
so it is natural to discard this information after some point.
The classical theory of packings retains only $a_{-1}$ (the density of $S$);
we suggest that it is natural to retain both $a_{-1}$ and $a_0$.
That is, we define a non-Archimedean valuation $\nu$
from the set of rational subsets of $\NN$ to $\QQ \times \QQ$,
where we view $\QQ \times \QQ$ as 
the lexicographic product of the ordered ring $\QQ$ with itself.
It can be shown that the pairs $(a_{-1},a_0)$ that occur
are those of the form $(0,k)$ or $(1,-k)$ where $k$ is a nonnegative integer,
along with pairs of the form $(p,q)$
where $p$ is a rational number strictly between 0 and 1
and where $q$ is an arbitrary rational number.
This valuation is not translation-invariant;
if $\nu(S) = (p,q)$, then $\nu(S+1) = (p,q-p)$.
Note that under this valuation, the sets
$\{3,6,9,12,15,18\}$ and $\{1,3,6,9,15,18\}$
discussed at the end of section 2 have the same size.
The valuation is emphatically not countably additive,
as can for instance be seen by viewing
$\NN$ as a union of singleton sets.

One can try to extend this valuation to various classes of sets
that include but are not limited to the rational subsets of $\NN$.
One way to do this without directly invoking
the expansion of $\sum_{n \in S} q^n$ as a Laurent series in $1-q$
is to define a partial preorder on the power set of $\NN$
(the {\em lim inf preorder})
such that $S$ dominates $S'$ in the lim inf preorder
iff $\liminf_{q \rightarrow 1^-} (\sum_{n \in S} q^n - \sum_{n \in S'} q^n) \geq 0$.
This partial preordering, restricted to the rational sets,
coincides with the total preordering obtained
by factoring the germ-ordering through the valuation $\nu$.

An important rationale for truncating the germs
comes from considering the role played by
the choice of regularization scheme.
If one wanted to extend our theory from packings in $\NN$ to packings in $\ZZ$
(with a view toward eventually looking at packings in $\RR^d$),
a different regularization scheme would be required
(since for $S \subseteq \ZZ$, $\sum_{n \in S} q^n$ diverges for all $q$ in (0,1)
unless $S$ is bounded below).
Two natural choices are
the germ of $\sum_{n \in S} q^{|n|}$ as $q \rightarrow 1^-$ (``$L^1$-regularization'')
and the germ of $\sum_{n \in S} q^{n^2}$ as $q \rightarrow 1^-$ (``$L^2$-regularization'').
It can be shown that, for rational sets $S \subseteq \ZZ$
(defined in the natural way from the monoid structure of $\ZZ$)
the pair $(a_{-1},a_0)$ is the same for $L^2$-regularization and $L^1$-regularization,
while later coefficients $a_n$ are different in the two theories.
Indeed, the valuation $\nu$ we constructed earlier,
mapping the set of rational subsets of $\NN$ to $\QQ \times \QQ$, is quite robust;
most sensible regularization schemes give rise to $\nu$.
This is just a restatement of the fact that the Grandi series and its variants
have the same value under most sensible ways of summing divergent series.

\bigskip

\noindent
{\bf 6. Connection to sphere-packing}

\bigskip

In the case of packing $\NN$ with translates of
$B = \{0,1,2,\dots,k-1\}$,
there is an appreciable {\em efficiency gap}
between the best packing and all other packings:

\medskip

\noindent
{\bf Theorem 8:} For $k \geq 1$ and $D = \{1,2,\dots,k-1\}$,
if $S^*$ is the $D$-avoiding set $\{0,k,2k,3k,\dots\}$
and $S$ is any other $D$-avoiding set,
$S_q \preceq (S^*)_q - \frac{1}{k} + O(1-q)$.

\medskip

\noindent
{\bf Proof:} We focus on the case $k=2$ for clarity.
Let $S^* = \{0,2,4,\dots\}$
and let $S$ be some $\{1\}$-avoiding set other than $S^*$.
We can write $S$ as the disjoint union of two sets,
one of the form $\{0,2,\dots,2(m-1)\}$ (empty if $m=0$)
and one of the form $\{t_1,t_2,t_3,\dots\}$
(with $t_1 < t_2 < t_3 < \dots$) satisfying
$t_1 \geq 2m+1$, $t_2 \geq 2m+3$, $t_3 \geq 2m+5$, etc.
The germ of $S$ is dominated by the germ of
$\{0,2,\dots,2(m-1)\} \cup \{2m+1,2m+3,2m+5,\dots\}$;
but this germ is the same (up to $O(1-q)$)
as the germ of $\{1,3,5,\dots\}$,
which falls short of the germ of $\{0,2,4,\dots\}$
by $\frac12 + O(1-q)$.
The case $k>2$ is similar. 
$\blacksquare$

\medskip

Packing problems and distance-avoidance problems in $\NN$
were chosen as a testbed for ideas about
packing problems and distance-avoidance problems in $\RR^n$,
and more specifically, sphere-packing problems.
Note that the problem of packing spheres of radius 1 in $\RR^n$
is equivalent to the problem of packing points in $\RR^n$
so that no two are at distance less than 2
(the points are the centers of the spheres).
We will not pursue the topic of sphere-packing here,
but we will mention the conjectures that motivated this work.

\medskip

\noindent
{\bf Conjecture 9:} Let $S$ be a subset of $\RR^2$,
no two of whose points are at distance less than 2,
and let $S^*$ be the set of center-points in
a hexagonal close-packing of disks of radius 1 in $\RR^2$. 
Let
$$\delta(S) = \liminf_{s \rightarrow \infty} \ \left( 
  \sum_{(x,y) \in S^*} e^{-(x^2+y^2)/s^2} 
- \sum_{(x,y) \in S} e^{-(x^2+y^2)/s^2} \right).$$
Then either $S$ is related to $S^*$ by an isometry of $\RR^2$, 
in which case $\delta(S) = 0$,
or else $S$ is not related to $S^*$ by an isometry of $\RR^2$, 
in which case $\delta(S) > 0$.

\medskip

\noindent
{\bf Remark:} In private communication, Henry Cohn has shown that
when $S$ is related to $S^*$ by an isometry of $\RR^2$,
$\delta(S)$ is indeed 0.

\medskip

\noindent
{\bf Conjecture 10:} In Conjecture 9,
``$\delta(S) > 0$'' can be replaced by ``$\delta(S) \geq 1$'' in the conclusion.

\medskip

The dichotomy between $\delta(S) = 0$ and $\delta(S) \geq 1$ in Conjecture 10
might at first seem to contradict
the continuity of the summands as a function of the positions of the points;
if all the points move continuously,
won't the lim inf also change continuously?
The catch is that the lim inf can (and often does) diverge.
For instance, if one obtains $S$ from $S^*$
by translating a half-plane's worth of points by $\epsilon > 0$,
or dilating the configuration $S^*$ by a factor of $c > 1$,
then the lim inf diverges, no matter how close $\epsilon$ is to 0,
or how close $c$ is to 1.

Clearly the bound in Conjecture 10 cannot be improved,
since removing a single point from $S^*$ gives a set $S$
for which the lim inf is exactly 1. 


\bigskip

\noindent
{\sc Acknowledgments:} This work has benefited from conversations with
Tibor Beke, Ilya Chernykh, Henry Cohn, David Feldman, Boris Hasselblatt, 
Alex Iosevich, Sinai Robins, and Omer Tamuz.

\bigskip

\noindent
{\bf References}

\medskip

\noindent
[Be] Vieri Benci, Emanuele Bottazzi, and Maura di Nasso,
``Elementary Numerosity and Measures'',
{\it J.\ Logic and Anal.} {\bf 6} (2014).

\medskip

\noindent
[Bl] Andreas Blass, Mauro Di Nasso, Marco Forti,
``Quasi-selective ultrafilters and asymptotic numerosities'',
{\it Adv. in Math.\ }{\bf 231} (2012) 1462--1486; 
\href{http://arxiv.org/abs/1011.2089}{{\tt http://arxiv.org/abs/1011.2089}}.

\medskip

\noindent
[Bo] Lewis Bowen and Charles Radin,
``Densest Packing of Equal Spheres in Hyperbolic Space'',
{\it Discrete Comput.\ Geom.\ }{\bf 29} (2003), 23--39.

\medskip

\noindent
[Ch] Ilya Chernykh, ``Non-Trivial Extension of Real Numbers'',
available at
\href{http://vixra.org/abs/1701.0617}{\tt http://vixra.org/abs/1701.0617}.

\medskip

\noindent
[Co] Henry Cohn, ``A Conceptual Breakthrough in Sphere Packing'',
{\it Notices of the AMS}, Volume 64, No.\ 2 (February 2017), 102--115.

\medskip

\noindent
[Ka] Fred Katz, ``Sets and Their Sizes'', 
\href{https://arxiv.org/abs/math/0106100}{\tt https://arxiv.org/abs/math/0106100}.

\medskip

\noindent
[Ku] Greg Kuperberg, ``Notions of Denseness'', 
{\it Geom.\ Topol.\ }{\bf 4} (2000) 277--292.

\end{document}